\documentclass{amsart}

\usepackage{epsfig}
\usepackage{amsmath, amsthm}
\usepackage{amscd}
\usepackage{latexsym}
\usepackage{mathrsfs}
\usepackage{bm}
\usepackage{hyperref}
\usepackage{graphicx}
\usepackage[dvipsnames,usenames]{color}
\hypersetup{
   bookmarks=true,         
   unicode=false,          
   pdftoolbar=true,        
   pdfmenubar=true,        
   pdffitwindow=false,     
   pdfstartview={FitH},    
   pdftitle={My title},    
   pdfauthor={Author},     
   pdfsubject={Subject},   
   pdfcreator={Creator},   
   pdfproducer={Producer}, 
   pdfkeywords={keywords}, 
   pdfnewwindow=true,      
   colorlinks=true,       
   linkcolor=blue,          
   citecolor=blue,        
   filecolor=magenta,      
   urlcolor=MidnightBlue          
}

\theoremstyle{plain}
\theoremstyle{plain}
\newtheorem{theorem}{Theorem}[section]
\newtheorem{cor}[theorem]{Corollary}
\newtheorem{lem}[theorem]{Lemma}
\newtheorem{pro}[theorem]{Proposition}
\newtheorem{Def}[theorem]{Definition}
\newtheorem{rem}[theorem]{Remark}
\numberwithin{equation}{section}

\newcommand{\af}{almost-Fuchsian}

\newcommand{\cc}{convex core}

\newcommand{\cis}{conformal infinities}
\newcommand{\cs}{conformal structure}
\newcommand{\cmc}{constant mean curvature}
\newcommand{\cmcs}{constant mean curvature surface}

\newcommand{\gc}{Gaussian curvature}

\newcommand{\Hd}{Hausdorff dimension}
\newcommand{\hqd}{holomorphic quadratic differential}

\newcommand{\hym}{hyperbolic metric}
\newcommand{\htm}{hyperbolic three-manifold}

\newcommand{\kg}{Kleinian group}
\newcommand{\ls}{limit set}
\newcommand{\mc}{mean curvature}

\newcommand{\ms}{minimal surface}

\newcommand{\ps}{parallel surface}
\newcommand{\pc}{principal curvature}
\newcommand{\qf}{quasi-Fuchsian}
\newcommand{\RS}{Riemann surface}
\newcommand{\rv}{renormalized volume}
\newcommand{\sff}{second fundamental form}

\newcommand{\TS}{Teichm\"{u}ller space}
\newcommand{\tm}{three-manifold}

\newcommand{\WP}{Weil-Petersson}
\newcommand{\WPm}{Weil-Petersson metric}

\newcommand{\be}{\begin{equation}}
\newcommand{\ene}{\end{equation}}
\newcommand{\bpf}{\begin{proof}}
\newcommand{\epf}{\end{proof}}
\newcommand{\bl}{\begin{lem}}
\newcommand{\el}{\end{lem}}
\newcommand{\bd}{\begin{Def}}
\newcommand{\ed}{\end{Def}}
\newcommand{\ben}{\begin{enumerate}}
\newcommand{\een}{\end{enumerate}}
\newcommand{\bp}{\begin{proof}}
\newcommand{\ep}{\end{proof}}
\newcommand{\bpo}{\begin{pro}}
\newcommand{\epo}{\end{pro}}
\newcommand{\beq}{\begin{equation*}}
\newcommand{\eeq}{\end{equation*}}
\newcommand{\bear}{\begin{eqnarray*}}
\newcommand{\eear}{\end{eqnarray*}}
\newcommand{\bt}{\begin{theorem}}
\newcommand{\et}{\end{theorem}}

\newcommand{\gbar}{\bar{g}}

\newcommand{\nablabar}{{\overline{\nabla}}}

\newcommand{\C}{\mathbb{C}}
\newcommand{\R}{\mathbb{R}}
\newcommand{\I}{\mathbb{I}}
\renewcommand{\H}{\mathbb{H}}

\newcommand{\Fcal}{\mathcal{F}}
\newcommand{\Acal}{\mathcal{A}}

\newcommand{\QF}{\mathcal{Q\!F}}

\newcommand{\Nscr}{\mathscr{N}}

\DeclareMathOperator{\dist}{dist}%
\DeclareMathOperator{\hyp}{hyp}%
\DeclareMathOperator{\PSL}{PSL}%
\DeclareMathOperator{\id}{id}%
\DeclareMathOperator{\Vol}{Vol}%
\DeclareMathOperator{\Isom}{Isom}%

\numberwithin{equation}{section}

\allowdisplaybreaks

\def\XXint#1#2#3{{\setbox0=\hbox{$#1{#2#3}{\int}$}
    \vcenter{\hbox{$#2#3$}}\kern-.5\wd0}}

\makeatletter
\def\@citestyle{\m@th\upshape\mdseries}
\def\citeform#1{{\bfseries#1}}
\def\@cite#1#2{{%
  \@citestyle[\citeform{#1}\if@tempswa, #2\fi]}}
\@ifundefined{cite }{%
  \expandafter\let\csname cite \endcsname\cite
  \edef\cite{\@nx\protect\@xp\@nx\csname cite \endcsname}%
}{}
\makeatother

\begin{document}

\title{On Almost-Fuchsian Manifolds}

\author{Zheng Huang}
\address{Department of Mathematics,
The City University of New York,
Staten Island, NY 10314, USA.}
\email{zheng.huang@csi.cuny.edu}

\author{Biao Wang}
\address{Department of Mathematics, Wesleyan University, Middletown, CT 06459, USA}
\email{bwang@wesleyan.edu}

\date{April 11, 2011}



\begin{abstract}
Almost-Fuchsian manifold is a class of complete {\htm}s.
Such a three-manifold is a {\qf} manifold which contains a
closed incompressible {\ms} with {\pc}s everywhere in the
range of $(-1,1)$. In such a manifold, the {\ms} is
unique and embedded, hence one can parametrize these
{\htm}s by their {\ms}s. In this paper we obtain estimates
on several geometric and analytical quantities of an {\af}
manifold $M$ in terms of the data on the {\ms}. In particular,
we obtain an upper bound for the hyperbolic volume of
the convex core of $M$, and an upper bound on the
Hausdorff dimension of the {\ls} associated to $M$. We
also constructed a {\qf} manifold which admits more than one {\ms}, and it
does not admit a foliation of closed surfaces of {\cmc}.
\end{abstract}

\maketitle


\section {Introduction}
Quasi-Fuchsian manifold is an important class of complete
hyperbolic {\tm}s. In hyperbolic geometry, {\qf}
manifolds and their moduli space, the quasi-Fuchsian space
$\QF(S)$, have been objects of extensive study in recent
decades. In particular, incompressible surfaces of small
{\pc}s play an important role in hyperbolic geometry and low
dimensional topology (\cite{Rub05}). Analogs of these surfaces
and {\tm}s are also a center of study in anti-de-Sitter
geometry (\cite{KS07, Mes07}). In this paper, we mostly
consider a subspace of the {\qf} space: {\af} manifolds.
They form a subspace of the same complex dimension $6g-6$
in $\QF(S)$, where $g \ge 2$ is the genus of any closed
incompressible surface in the manifold (\cite {Uhl83}).
Understanding the structures of the quasi-Fuchsian space is a
mixture of understanding the geometry of quasi-Fuchsian
$3$-manifolds, the deformation of incompressible surfaces,
as well as the representation theory of {\kg}s. It is highly
desirable to use information on special surfaces (minimal or
{\cmc}) to obtain global information on the {\tm}.

Let us make a definition:
\vskip 0.1in
{\bf Definition:} {\it
We call $M$ an {\af} manifold if it is a {\qf} manifold which
contains a closed incompressible {\ms} $\Sigma$ such that
the {\pc}s of $\Sigma$ are in the range of $(-1,1)$}.
\vskip 0.in
Recall that we call a closed surface {\it incompressible} in
$M$ if the inclusion induces an injection between fundamental
groups of the surface and the {\tm} $M$. Naturally, $M$ is
Fuchsian when $\Sigma$ is actually totally geodesic. The
notion of {\af} (term coined in \cite{KS07}) was first studied
by Uhlenbeck (\cite{Uhl83}), where she proved several key
properties of {\af} manifolds that will be vital in this work:
$\Sigma$ is the only closed incompressible {\ms} in $M$, and
$M$ admits a foliation of {\ps}s from $\Sigma$ to both ends.
Around the same time, Epstein (\cite{Eps86}) used the
parallel flow in $\H^3$ to study the quasiconformal reflection
problem.

Throughout the paper, all surfaces in $M$ involved are assumed to
be closed, oriented, incompressible, of genus at
least two. We also assume $M$ is not Fuchsian, or most theorems
are trivial. One of our motivations is to use various
techniques in analysis to investigate geometric problems in {\htm}s.

Since the {\ms} is unique in an {\af} manifold, we can use {\ms}s
to parametrize the space of {\af} manifolds. This
parametrization is in terms of the {\cs} of the {\ms} and the
{\hqd} on the {\cs} which determines the {\sff} of the {\ms}
in the {\tm} (see \cite{Uhl83, Tau04, HL12}). We will focus on
instead obtaining topological and geometric information
about $M$ from data of $\Sigma$ in this paper. Among the quantities
associated to a {\qf} manifold, the hyperbolic
volume of the {\cc} and the {\Hd} of the limit set are perhaps the
most significant. Estimates on them in term
of the geometry of the {\ms} in $M$ are obtained in this paper.
Moreover, the foliation structure possessed by $M$
allows us to investigate a notion of renormalized volume for the
{\af} manifold.

The {\cc} of a {\qf} manifold is the smallest convex subset of a
{\qf} manifold that carries its fundamental group. From
the point of view of hyperbolic geometry, the {\cc} contains all
the geometrical information about the {\qf} {\tm} itself
(see for instance, \cite {AC96, Bro03}). As a direct application,
when $M$ is {\af}, we obtain an explicit upper bound
for the hyperbolic volume of the {\cc} $C(M)$, in terms of the maximum
{\pc} on the {\ms} $\Sigma$:

\bt\label{thm:volume estimate}
If $M$ is almost-Fuchsian, and let
$\lambda_0 = \max\limits_{x \in \Sigma}\{|\lambda(x)|\}$ be the maximum of the
{\pc} of the {\ms} $\Sigma$ of $M$, then
\begin{align}
    \Vol(C(M))
       &\leq \Acal_{\hyp}\left({\frac{\lambda_0}
            {1-\lambda_0^2}}+{\frac{1}{2}}
            \log{\frac{1+\lambda_0}
            {1-\lambda_0}}\right) \notag \\
       &= \Acal_{\hyp} \left(2\lambda_0 +
          {\frac{4}{3}}\lambda_0^3 +
          O(\lambda_0^5)\right),
\end{align}
where $\Acal_{\hyp} = 2\pi(2g-2)$ is the hyperbolic area of $\Sigma$.
\et

We immediately see quantitatively how the volume of $C(M)$ goes to zero
when $\lambda_0$ is close to zero. In
\cite{Bro03}, Brock showed the hyperbolic volume of the {\cc} is
quasi-isometric to the {\WP} distance between {\cis}
of $M$ in {\TS}.

It is well-known that the {\Hd} of the limit set for any {\qf} group
is in the range of $[1,2)$, and identically $1$ if and only
if it is Fuchsian. We also obtain an upper bound for the {\Hd} of the
limit set $\Lambda_{\Gamma}$ of $M$, in terms of
$\lambda_0$ as well:

\bt\label{thm:dimension estimate}
If $M= \H^3/\Gamma$ is {\af}, and let
$\lambda_0 = \max\limits_{x \in \Sigma}\{|\lambda(x)|\}$, then the
{\Hd} $D(\Lambda_{\Gamma})$
of the limit set $\Lambda_{\Gamma}$ for $M$ satisfies
\be
   D(\Lambda_{\Gamma}) < 1+ \lambda_0^2.
\ene
\et

For $\lambda_0$ close to zero, Theorems \ref{thm:volume estimate} and
\ref{thm:dimension estimate} measure how
close $M$ is to being Fuchsian. In (\cite{GHW10}), we related the
foliation structure of an {\af} manifold to both
Teichm\"uller metric and the {\WPm} on {\TS}.

Given a {\qf} manifold, it admits finitely number of {\ms}s (\cite{And83}).
One expects the space of {\af} manifolds is not
the full {\qf} space, and indeed, the second named author (\cite{Wan12})
showed some {\qf} manifolds that admit
more than one {\ms}. A further generalization is to consider closed
surfaces of {\cmc} in a {\qf} manifold. This is a
vastly rich area where many analytical techniques can be applied.
A natural question (Thurston) is to ask, to what
extent that a {\qf} manifold admits a foliation of closed (incompressible)
surfaces of {\cmc}. In the second part of this
paper, we construct a {\qf} manifold which admits more than one {\ms} (hence not {\af}),
such that it does not admit such a foliation, namely,

\bt\label{non}
There exists a {\qf} manifold $N$ which does not admit a foliation of {\cmcs}s.
\et

\vskip 0.1in
One of our original motivations is to investigate whether {\af}
manifolds are the appropriate subclass of {\qf} manifolds
that admit such a foliation.

We conclude this introduction by the following note: our volume estimate
can be generalized to a possibly slightly larger class
of {\qf} manifolds than {\af} manifolds. In other words, one can
define a notion of {\it nearly Fuchsian} manifolds: a class
of {\qf} manifolds that each admits a closed incompressible surface
(not necessarily minimal) of {\pc}s in the range of
$(-1,1)$. One can similarly verify that a nearly Fuchsian manifold admits a foliation
by {\ps}s of this fixed surface of small {\pc}s. It is not known if these two classes of {\qf} manifolds
actually coincide, or if any nearly Fuchsian manifold admit only one {\ms}.

\subsection*{Plan of the paper}
After a brief section on the preliminaries, the rest of the paper
contains two parts: In \S \ref{sec:AF manifolds geometry},
we prove several results on the geometry of {\af} manifolds,
including the {\rv} of $M$, Theorem
\ref{thm:volume estimate} (volume estimate) and Theorem
\ref{thm:dimension estimate} ({\Hd} estimate). In \S
\ref{sec:counterexam}, we prove Theorem \ref{non} by
constructing an explicit {\qf} manifold.

\subsection*{Acknowledgements}
The authors are grateful to Ren Guo for many helpful discussions,
Dick Canary and Jun Hu for their suggestions regarding
the {\Hd} of the limit set. The last section is based on part of
the thesis of the second named author, and he wishes to thank
his advisor Bill Thurston for advices, inspiration and encouragement.
We also want to express our gratitude to anonymous referee(s), for their
willingness to read the paper carefully and for many suggestions to improve this paper.
The research of Zheng Huang is partially supported by a PSC-CUNY
research award and an award from the CUNY-CIRG program and CUNY-CSI Provost's Scholarship.


\section{Preliminaries}\label{sec:prelim}
In this section, we fix our notations, and introduce some
preliminary facts that will be used later in this paper.


\subsection{Quasi-Fuchsian manifolds}
For detailed reference on Kleinian groups and low dimensional
topology, one can go to \cite{Mar74} and \cite{Thu82}.

The universal cover of a complete orientable {\htm} is $\H^3$, and the deck transformations
induce a representation of the fundamental
group of the manifold in $\Isom(\H^3)=\PSL(2,\C)$, the (orientation
preserving) isometry group of $\H^3$. A subgroup
$\Gamma \subset \PSL(2,\C)$ is called a {\em Kleinian group} if
$\Gamma$ acts on $\H^{3}$ properly discontinuously.
For any {\kg} $\Gamma$, $\forall\,p\in\H^{3}$, the orbit set
\begin{eqnarray*}
   \Gamma(p)=\{\gamma(p)\ |\ \gamma\in \Gamma\}
\end{eqnarray*}
has accumulation points on the boundary
$\partial\H^{3}=S^{2}_\infty=\widehat{\C}$,
and these points are the {\em limit points}
of $\Gamma$, and the closed set of all these points is called the
{\em limit set} of $\Gamma$, denoted by
$\Lambda_{\Gamma}$.

In the case when $\Lambda_{\Gamma}$ is contained in a circle
$S^1 \subset S^2$, the quotient manifold
$M =\H^{3}/\Gamma$ is called {\it Fuchsian}, and $M$ is isometric
to a (warped) product space $S \times \R$.
If the limit set $\Lambda_{\Gamma}$ lies in a Jordan curve, the
quotient {\tm} $M =\H^{3}/\Gamma$ is called {\em {\qf}},
and it is topologically $S \times \R$, where $S$ is a closed surface
of genus $g$ at least two. It is clear that a {\qf}
manifold is quasi-isometric to a Fuchsian manifold. The space of
such {\tm}s $\QF(S)$, the {\qf} space of genus $g$ surfaces,
is a complex manifold of dimension of $6g-6$, which has very
complicated structures.

Finding {\ms}s in negatively curved manifolds is a problem of
fundamental importance. The basic results are due to
Schoen-Yau (\cite {SY79}) and Sacks-Uhlenbeck (\cite{SU82}),
and their results can be applied to the
case of {\qf} manifolds: any quasi-Fuchsian manifold contains
at least one incompressible {\ms}. In the case of
almost-Fuchsian, the {\ms} is unique (\cite{Uhl83}).
On the other hand, there are {\qf} manifolds that
admit many minimal surfaces (\cite{Wan12}).

An essential problem in hyperbolic geometry and complex dynamics
is to study the {\Hd} $D(\Lambda_{\Gamma})$ of
the limit set $\Lambda_{\Gamma}$ associated to $M$. This problem
is also intimately related to understanding the
lower spectrum theory of {\htm} (\cite{Sul87, BC94}). In the case
of Fuchsian manifolds, $\Lambda_{\Gamma}$ is a
round circle and $D(\Lambda_{\Gamma}) = 1$. When $M$ is {\qf} but
not Fuchsian, as is throughout this paper, it is
known that $1 < D(\Lambda_{\Gamma}) < 2$ (\cite{Bow79, Sul87}). There is a rich theory
of quasiconformal mapping and its distortion on
{\Hd}, area and other quantities (see for instance \cite {GV73, LV73}).


\subsection{Almost-Fuchsian manifolds}
We now assume that $M$ is an {\af} manifold: the {\pc}s of the {\ms}
$\Sigma$ are in the range of $(-1,1)$. It is clear
that for any closed embedded surface $S$ in $M$, one can define
a regular parallel surface $S(r)$ which is fixed
(hyperbolic) distant $r$ from $S$, for sufficiently small $r$.
A remarkable property for $M$ is that when taking $S$
to be the {\ms}, the parallel surface $S(r)$ is nonsingular for
all $r \in \R$.

To be more precise, using isothermal coordinates, the induced metric on a closed incompressible
surface $S$ is given by
$g_{ij}(x)=e^{2v(x)}\delta_{ij}$, where $v(x)$ is a smooth
function on $S$, and while the {\sff} is denoted by
$A(x)=[h_{ij}]_{2\times{}2}$, where $h_{ij}$ is given by, for
$1\leq{}i,j\leq{}2$,
\begin{equation*}
   h_{ij}=\langle{\nablabar_{e_{i}}\nu}, {e_{j}}\rangle,
\end{equation*}
where we choose $\{e_{1},e_{2}\}$ as an orthonormal basis on $S$, and $\nu$
is the unit normal field on $S$, and $\nablabar$ is the Levi-Civita
connection of $(M,\gbar_{\alpha\beta})$. Here, we add a bar on top
for each quantity or operator with respect to $(M,\gbar_{\alpha\beta})$.

Let $\lambda_1(x)$ and $\lambda_2(x)$ be the eigenvalues of $A(x)$.
They are the {\pc}s of $S$, and we denote
$H(x) =\lambda_1(x)+ \lambda_2(x)$ as the {\mc} function of $S$.
In classical differential geometry, the {\sff} indicates
how a hypersurface immerses into the ambient manifold. Zero {\sff}
is equivalent to say that the hypersuface is totally
geodesic. Therefore, {\pc}s are natural quantities to investigate
in this type of problems.

Let $S(r)$ be the family of equidistant surfaces with respect
to $S$, i.e.
\begin{equation*}
   S(r)=\{\exp_{x}(r\nu)\ |\ x\in{}S\}\ ,
   \quad{}r\in(-\varepsilon,\varepsilon)\ .
\end{equation*}
The induced metric on $S(r)$ is denoted by
$g(x,r) = [g_{ij}(x,r)]_{1\leq{}i,j\leq{}2}$, and the {\sff} is
denoted by $A(x,r)=[h_{ij}(x,r)]_{1\leq{}i,j\leq{}2}$. The {\mc}
on $S(r)$ is thus given by $H(x,r)=g^{ij}(x,r)h_{ij}(x,r)$. Uhlenbeck
calculated the induced metric on $S(r)$ in terms of the induced metric
on $S$ and the {\sff} of $S$, and found

\bl[\cite {Uhl83}]
The induced metric $g(x,r)$ on $S(r)$ has the form
\be\label{eq:metric of S(r)}
   g(x,r)=e^{2v(x)}[\cosh(r)\I+\sinh(r)
   e^{-2v(x)}A(x)]^{2}\ ,
\ene
where $r \in (-\varepsilon,\varepsilon)$.
\el

Any point on the parallel surface $S(r)$ in $M$ can now be
represented by a pair $(x,r)$, where $x \in S$ and
$r \in (-\varepsilon,\varepsilon)$. Direct computation shows
that the {\pc}s of $S(r)$ are given by
\begin{equation}\label{eq:pc4ef}
   \mu_{j}(x,r)=\frac{\tanh{}r+\lambda_{j}(x)}
   {1+\lambda_{j}(x)\tanh{}r}\ ,\qquad j=1,2\ .
\end{equation}
Hence the {\mc} (the sum of {\pc}s) is given by
\be \label{eq:mc4ef}
   H(x,r)=\frac{2(1+\lambda_{1}\lambda_{2})\tanh{}r+
   (\lambda_{1}+\lambda_{2})(1+\tanh^{2}r)}
   {1+(\lambda_{1}+\lambda_{2})\tanh{}r+\lambda_{1}
   \lambda_{2}\tanh^{2}r}\ .
\ene

Since $M$ is {\af}, we now choose $S = \Sigma$ to be
the unique {\ms} in $M$. From the explicit nature of the
Lemma \ref{eq:metric of S(r)}, one concludes that, when
$|\lambda_{j}(x)| < 1$ for $j=1,2$ and $x \in \Sigma$,
the induced metrics $g(x,r)$ are of no singularity for
all $r \in \R$ and therefore {\ps}s of $\Sigma$ form a
foliation on $M$, called the {\it equidistant foliation}
or the {\it normal flow}. We also observe all {\pc}s
$\{\mu_j(x,r)\}$ on the {\ps}s are in the range of $(-1,1)$.
We denote the equidistant foliation from the {\ms}
$\Sigma$ by $\{\Sigma(r)\}_{r\in \R}$.


\section{Geometry of Almost-Fuchsian Manifolds}
\label{sec:AF manifolds geometry}

In this section, we wish to obtain information about the {\af}
manifold $M$ via its unique {\ms} $\Sigma$. We
derive several geometrical properties on the equidistant foliation
$\{\Sigma(r)\}_{r \in \R}$ in \S
\ref{subsec:some estimate}. The calculations in this subsection
will set up the stage for later applications.
In \S \ref{subsec:upper bound}, we obtain explicit upper bounds
for the hyperbolic volume of the {\cc} of $M$, and compute explicitly the {\rv} of $M$. In \S 3.3,
we establish the estimate for the {\Hd} of the {\ls} associated to the {\af} manifold $M$.

\subsection{Some estimates on the {\ps}s}
\label{subsec:some estimate}
In this subsection, we provide several estimates that will be used
later. First, let us record a few quantities
that will be involved:
\ben
\item The {\pc}s of the minimal surface $\Sigma$ are
      $\pm \lambda(x)$, where $x \in \Sigma$ and we have
      $0 \le \lambda(x) < 1$;
\item $\lambda_0$ is the maximum of $\lambda$, hence the
      maximal {\pc} on $\Sigma$;
\item $|S|$ is the area for any closed incompressible surface
      $S$ (with respect to the induced metric);
\item $\Acal_{\hyp} = 2\pi(2g-2)$ is the hyperbolic area
      of any closed {\RS} $S$;
\item $J(r)$ is the determinant of the Jacobian between
      $\Sigma$ and the {\ps} $\Sigma(r)$ via the pullback.
      Note that $\Sigma(r)$ is regular if and only if $J(r) > 0$;
\item $K(r)$ is the {\gc} on $\Sigma(r)$. In particular,
      $K(0) = K(\Sigma)$ is the {\gc} on the {\ms} $\Sigma$;
\item $M(\pm r)$, $r > 0$, is the region in $M$ bounded by
      surfaces $\Sigma(-r)$ and $\Sigma(r)$;
\item Lastly, $M(r)$, $r > 0$, is the region in $M$ bounded
      between the {\ms} $\Sigma$ and the {\ps} $\Sigma(r)$.
\een

We start with a well-known estimate which implies the area of the
{\ms} under the induced metric from the
ambient space is comparable to that of the hyperbolic area, with
universal constants. We only include a proof
for the sake of completeness, and it is very short.

\bpo\label{prop:vol estimate 1}
$\frac{\Acal_{\hyp}}{2} < |\Sigma| < \Acal_{\hyp}$.
\epo

\begin{proof}
We apply the Gauss equation:
\begin{equation*}
   K(\Sigma) = -1 + \det(A) = -1 -\lambda^2
\end{equation*}
Thus we have
\begin{equation*}
   -K(\Sigma) = 1 - \det(A) = 1 + \lambda^2.
\end{equation*}
We integrate this on the surface $\Sigma$, and apply
the Gauss-Bonnet theorem, since $\Sigma$ is closed, to find
\begin{equation*}
   |\Sigma| < |\Sigma| + \int_{\Sigma}\lambda^2 =
   \Acal_{\hyp} < 2|\Sigma|.
\end{equation*}
\end{proof}

We next want to estimate the area of each {\ps} in the equidistant
foliation $\{\Sigma(r)\}_{r \in \R}$. We can see
that they grow at a rate of $\sinh^2(r)$, which is as expected.
The explicit formula in Lemma \ref{eq:metric of S(r)}
indicates, for large $|r|$, the metric for $M$ behaves like the
warped product metric.

\bpo\label{prop:vol estimate 2}
For all $-\infty < r < +\infty$, we have
\begin{equation*}
  (2|\Sigma|-\Acal_{\hyp})\sinh^2(r) <
  |\Sigma(r)| < |\Sigma |\cosh^2{}r < \Acal_{\hyp}\cosh^2{}r
\end{equation*}
\epo

\begin{proof}
The area element of $\Sigma (r)$ is given by
\be\label{jacobian}
d\mu(r) = (\cosh^2{}r - \lambda^{2}(x)\sinh^{2}{}r)d\mu,
\ene
where $d\mu$ is the area element for the {\ms} $\Sigma$.

We can now compute the surface area:
\begin{align}
  |\Sigma(r)| &= \int_{\Sigma} (\cosh^2{}r -
            \lambda^{2}(x)\sinh^{2}{}r)d\mu  \nonumber \\
         &= |\Sigma|\cosh^{2}{}r  - \sinh^{2}{}r
             \int_{\Sigma}\lambda^{2}(x)d\mu \nonumber \\
             \label{eq:surface area (r)-1}
         &= |\Sigma|\cosh^{2}{}r  -
            (\Acal_{\hyp} - |\Sigma|)\sinh^{2}{}r \\
         &= |\Sigma|(\cosh^{2}{}r + \sinh^{2}{}r) -
            \Acal_{\hyp}\sinh^{2}{}r\nonumber \\
            \label{eq:surface area (r)-2}
         &= |\Sigma| + (2|\Sigma|-\Acal_{\hyp})\sinh^{2}{}r \,
\end{align}
Here we used the identity
\be
\int_{\Sigma}\lambda^2 = \Acal_{\hyp} - |\Sigma|.
\ene
The estimates then follows from the Proposition
\ref{prop:vol estimate 1}, formulas
\eqref{eq:surface area (r)-1} and \eqref{eq:surface area (r)-2}.
\end{proof}

As a consequence of the formula \eqref{jacobian}, we find that the determinant of the
Jacobian between $\Sigma$ and $\Sigma(r)$ is
\be
J(r) = \cosh^2(r) - \lambda^{2}(x)\sinh^{2}(r).
\ene
From this, clearly, $|\lambda| < 1$ implies that $J(r) > 0$
for all $r \in \R$. We also have the following
version of the Gauss-Bonnet formula:

\bpo
The product of $J(r)$ and $K(r)$ is a function of $x \in \Sigma$,
independent of $r$. In other words,
\be\label{JK}
J(r)K(r) = K(0) = -1-\lambda^2.
\ene
\epo
\bp

Let us conduct this computation. The formulas for the
{\pc}s of $\Sigma(r)$ are given by \eqref{eq:pc4ef}.
Therefore we have
\bear
K(r) & = & -1 + \mu_1(r)\mu_2(r)\\
& = & -1 + \frac{\tanh^2(r)-\lambda^2}{1-\lambda^2 \tanh^2(r)} \\
& = & \frac{(\lambda^2+1)(\sinh^2(r) - \cosh^2(r))}{J(r)} \\
& = & -\frac{(\lambda^2+1)}{J(r)}.
\eear
\ep
We conclude this subsection with the following remark that the
identity \eqref{JK} holds more generally. If
$S$ is a closed surface with {\pc}s $\{\lambda_1(x),\lambda_2(x)\}$,
and $r \in \R$ is a real number such
that $S(r)$ is a nonsingular {\ps} of $S$. Then the identity
\eqref{JK} shows that what must happen
geometrically when $S(r)$ is developing a singularity:
the {\gc} of $S(r)$ must blow up. This was also observed in \cite{Eps84}.

\subsection{The {\cc} volume and the {\rv}}
\label{subsec:upper bound}

We obtain an upper bound for the hyperbolic volume of the {\cc} $C(M)$ in this subsection, in terms of the
maximum, $\lambda_0$, of the {\pc} on the {\ms} $\Sigma$. The idea and actual computation are
somewhat simple: we take advantage of the foliation structure of the {\af} manifold $M$, and note that, from
the formula \eqref{eq:pc4ef}, the {\pc}s of the surface $\Sigma(r)$, $\mu_{1}(x,r)$ and $\mu_{2}(x,r)$ are
increasing functions of $r$ for any fixed $x \in \Sigma$, and they approach $\pm 2$ as
$r \rightarrow \pm\infty$. When any of the {\ps}s becomes convex (all positive {\pc}s or all negative {\pc}s),
they lie outside of the {\cc}.

Naturally, we are particularly interested in two critical cases: the values of $r$ when $\mu_1(x,r) = 0$ or
$\mu_2(x,r) = 0$. Elementary algebra shows:

\bpo\label{r_0}
 If we denote
\begin{equation*}
   r_0 = {\frac{1}{2}}\log{\frac{1+\lambda_0}{1-\lambda_0}}\ ,
\end{equation*}
where $\lambda_0 = \max\limits_{x \in S}\{\lambda(x)\}$, then $r_0$ is the least value of $r$ such that
$\mu_1(r,x) > 0$ for all $r > r_0$ and $x \in \Sigma$, while $-r_0$ is the largest value for $\mu_2(r,x) < 0$
such that $\mu_2(r,x)<0$ for all $r<-r_0$ and $x \in \Sigma$.
\epo

This Proposition tells us when the {\ps}s in the equidistant foliation $\{\Sigma(r)\}_{r \in \R}$ become
convex, hence by the definition of the convex core, provides an upper bound for the size of the {\cc}.

Recalling that we denote the region of $M$ bounded between surfaces $\Sigma(-r_0)$ and $\Sigma(r_0)$
by $M(\pm r_0)$, and then the convex core $C(M)$, is contained in $M(\pm r_0)$. Since
$\{\Sigma(r)\}_{r \in \R}$ foliates $M$, we can compute the hyperbolic volume of the region $M(\pm r_0)$
by the following:
\begin{align}
   \Vol(M(\pm r_0))
       &= \int_{-r_0}^{r_0}|\Sigma(r)|dr \nonumber \\
       &= 2r_{0}|\Sigma| + (2|\Sigma| - \Acal_{\hyp})
          \int_{-r_0}^{r_0}\sinh^{2}r{}dr \nonumber \\
       &= 2r_{0}|\Sigma| + (2|\Sigma| - \Acal_{\hyp})
          \left({\frac{1}{2}}\sinh(2r_0) -
          r_0\right)  \nonumber \\
          \label{eq:vol estimate r0}
       &= |\Sigma|\sinh(2r_0) - \Acal_{\hyp}
          \left({\frac{1}{2}}\sinh(2r_0) - r_0\right)
\end{align}
Applying the Proposition \ref{prop:vol estimate 2}, we obtain the following:

\bt\label{thm:hyperbolic volume estimate}
The hyperbolic volume of $C(M)$ is bounded by:
\begin{eqnarray}
    \Vol(C(M)) &\le& \Acal_{\hyp}(\cosh{}r_{0}
                    \sinh{}r_{0} + r_0)\nonumber \\
               &=& \Acal_{\hyp}\left({\frac{\lambda_0}
                  {1-\lambda_0^2}}+
                  {\frac{1}{2}}\log{\frac{1+\lambda_0}
                  {1-\lambda_0}}\right).
\end{eqnarray}
\et
\medskip
When $r_0 = 0$, or equivalently, $\lambda(x) =0$ for all $x \in \Sigma$, this is the case of $M$ being
Fuchsian, and the hyperbolic volume of the {\cc} $C(M)$ is zero. We want to measure how the
hyperbolic volumes vary for small $\lambda_0$ via the following Taylor series expansion:

\begin{cor}
For small $\lambda_0$, we have the following expansion:
\begin{equation*}
   \Vol(C(M)) \le \Acal_{\hyp}\left(2\lambda_0 +
   {\frac{4}{3}}\lambda_{0}^{3} + O(\lambda_0^5)\right).
\end{equation*}
\end{cor}

Any {\qf} manifold is complete, hence has infinite volume.
It is classical in conformal geometry to define
a notion of {\rv} (\cite{FG85, PP01}) to obtain some conformal invariant.
The foliation structure of an
{\af} manifold allows one to derive a simple quantity as the {\rv}.
We adapt the following notion: the
{\it \rv} of $M$ with respect to the foliation $\{\Sigma(r)\}_{r \in \R}$ is given by
\be\label{lim}
RV(M) = RV(M, \{\Sigma(r)\}) =
2\lim_{r\to\infty}\left\{\Vol(M(r)) - \frac12 |\Sigma(r)| - \frac{\Acal_{\hyp}}{2}r\right\},
\ene
where we recall that $M(r)$ is the region bounded by $\Sigma$ and $\Sigma(r)$. Note here any {\qf}
manifold has two ends and we take advantage of our situation of the obvious symmetry of the foliation
$\{\Sigma(r)\}_{r \in \R}$ with respect to the {\ms} $\Sigma$. One of the applications from computing
the volume of $M(r)$ is to determine this limit:
\bpo
When $M$ is {\af}, the {\rv} (with respect to the foliation $\{\Sigma(r)\}_{r \in \R}$) is
\beq
RV(M) = 2\pi (1-g).
\eeq
\epo
\bpf
We will collect and organize terms in the limit. Firstly, using \eqref{eq:vol estimate r0}, we have the following:
\beq
\Vol(M(r)) = \frac12 \Vol(M(\pm r)) =
\frac14 \sinh(2r)(2|\Sigma| - \Acal_{\hyp}) + \frac{\Acal_{\hyp}}{2}r.
\eeq
Secondly, we apply \eqref{eq:surface area (r)-2} to find
\beq
 \frac12 |\Sigma(r)|= \frac12 |\Sigma| + (2|\Sigma|-\Acal_{\hyp})\frac{\sinh^{2}(r)}{2}.
\eeq
Combining these terms, we have
\bear
\Vol(M(r)) - \frac12 |\Sigma(r)| - \frac{\Acal_{\hyp}}{2}r & = & (2|\Sigma|-\Acal_{\hyp})\left(\frac{1-e^{-2r}}{4}\right)-\frac12 |\Sigma|\\
&=& -\frac{\Acal_{\hyp}}{4} + \frac{2|\Sigma|-\Acal_{\hyp}}{4}e^{-2r}.
\eear
Now the statement holds after taking the limit for $r \to \infty$ and using \eqref{lim}.
\epf

This limit, as a volume, may sound unnatural for it is negative in this case. This is however typically the case when one
try to get a finite quantity out of a diverging sequence: one expands the (unbounded) quantity with respect to a parameter,
and obtain the desired finite quantity from the constant term in the expansion. In our case, the {\rv} is the constant term in
the series expansion of $\Vol(M(r))$ in terms of $r$ for large $r > 0$. One can regard this as the mass is negative for the {\af} manifolds.

Also note that $H(r) \to 2$ as $r \to \infty$, in other words, the {\ps}s $\Sigma(r)$ are close to {\cmc} surfaces as $r$ gets large.
Therefore one might attempt to use the integral $\int_{\Sigma(r)}H(r)$ to replace the term $|\Sigma(r)|$ in
\eqref{lim}. This offers an alternative interpretation for the {\rv} of $M$: it is the limit of half of the total
difference of the {\mc} $H(r)$ and $2$ (the {\cmc} of the infinity), on $\Sigma(r)$, for large $r$:
\begin{cor}(also \cite{KS08})
\be
RV(M) = \frac12\lim_{r\to\infty}\left\{\int_{\Sigma(r)}(H(r) -2)\right\}.
\ene
\end{cor}
\bpf
This follows from the following identity: for $r > 0$,
\be\label{vol}
\Vol(M(r)) = \frac14\int_{\Sigma(r)} H(r) + \frac{r}{2}\Acal_{\hyp}.
\ene
This identity can be easily verified by the following:
\bear
\int_{\Sigma(r)} H(r) &=& \int_{\Sigma}(\mu_1(r) + \mu_2(r))J(r)d\mu\\
&=&\sinh(2r)(2|\Sigma| - \Acal_{\hyp}).
\eear
\epf
More generally, the {\rv} for {\qf} manifolds was investigated in \cite{KS08}, where they obtained more general
identities than \eqref{vol} relating the {\rv} and the total {\mc} of {\ps}s.
\subsection{Hausdorff dimension of the limit set} A {\qf} manifold is determined by a subgroup $\Gamma$
of $\PSL(2,\C)$. It is a natural question to ask how much one knows about the group $\Gamma$ when the
resulting {\qf} manifold $M = \H^3/\Gamma$ is {\af}. In this subsection, we attempt to investigate this
connection. A critical question is to understand the limit set of $\Gamma$, and in our case, to determine
its {\Hd} via a geometric quantity. We obtain two estimates. One is a straightforward application of our
prior volume estimate for the {\cc} and a theorem of Burger-Canary (\cite{BC94}), and the other approach
is more technical, but with a much simpler answer.

We proceed with the first approach. We denote $C_1(M)$ the hyperbolic radius one neighborhood of the
{\cc} $C(M)$ in $M$. An easy calculation from \eqref{eq:vol estimate r0} and Proposition
\ref{prop:vol estimate 1} show us:
\beq
   \Vol(C_1(M)) \leq 2\Vol(M(r_0+1)) \leq  \Acal_{\hyp}\left({\frac{1}{2}}\sinh(2r_0+2) + r_0 +1\right),
\eeq
where  $r_0 = {\frac{1}{2}}\log{\frac{1+\lambda_0}{1-\lambda_0}}$. Therefore we have
\begin{equation}
   \Vol(C_1(M)) \le  \Acal_{\hyp}\left({\frac{1}{2}}
   \sinh\left(\log{\frac{1+\lambda_0}{1-\lambda_0}}+2\right) +
   {\frac{1}{2}}\log{\frac{1+\lambda_0}{1-\lambda_0}} +1\right).
\end{equation}
Since quasi-Fuchsian manifolds are geometrically finite and of infinite volume, and we assume $M$ is
not Fuchsian, a direct application of the main theorem from Burger-Canary (\cite{BC94}) gives:

\bpo
Let $M$ be {\af}, and $\mu_0(M)$ be the bottom of the $L^2$-spectrum of $-\Delta$ on $M$, and
$D(\Lambda_{\Gamma})$ be the {\Hd} of the limit set $\Lambda_{\Gamma}$ of $M$. Then we have
\begin{enumerate}
\item
\begin{equation}
   \mu_0(M) \ge  {\frac{K_3}{\Acal_{\hyp}^2\left({\frac{1}{2}}
   \sinh(\log{\frac{1+\lambda_0}{1-\lambda_0}}+2) +
   {\frac{1}{2}}\log{\frac{1+\lambda_0}{1-\lambda_0}} +1\right)^2}}.
\end{equation}
\item
\begin{equation}
   D(\Lambda_{\Gamma}) \le  2-{\frac{K_3}{\Acal_{\hyp}^2\left({\frac{1}{2}}
   \sinh(\log{\frac{1+\lambda_0}{1-\lambda_0}}+2) +
   {\frac{1}{2}}\log{\frac{1+\lambda_0}{1-\lambda_0}} +1\right)^2}}.
\end{equation}
\end{enumerate}
Here $K_3$ can be chosen such that $K_3 > 10^{-11}$.
\epo
\medskip
We note that while the volume estimate of the {\cc} of $M$ (Theorem \ref{thm:hyperbolic volume estimate})
is effective for small maximal {\pc} $\lambda_0$ of the {\ms} $\Sigma$, above estimates on $\mu_0(M)$
and $D(\Lambda_{\Gamma})$ are not as effective. To obtain an estimate only depending on the {\ms}, we
switch to a different approach: we consider the limit set $\Lambda_{\Gamma}$ of $M$ as a $k$-quasicircle
(an image of a circle under a $k$-quasiconformal mapping), and our task is reduced to estimate $k$ in terms of
$\lambda_0$.

A {\it $k$-quasiconformal} mapping $f$ is a homemorphism of planar domains, locally in
the Sobolev class $W_2^1$ such that its {\it Beltrami coefficient}
$\mu_f = \frac{\bar\partial f}{\partial f}$ has bounded $L^{\infty}$ bound:
$\|\mu_f\| \le k < 1$. One can visualize $f$ infinitesimally maps a round circle to an ellipse
with a bounded dilatation $K = \frac{1+k}{1-k}$, where $k \in [0,1)$. Clearly, the mapping
$f$ is conformal when $k =0$.

We now prove the Theorem \ref{thm:dimension estimate}, which we re-state here:
\bt
Let $M$ be {\af}, then the {\Hd} $D(\Lambda_{\Gamma})$ of the limit set $\Lambda_{\Gamma}$ for
$M = \H^3/\Gamma$ satisfies
\beq
   D(\Lambda_{\Gamma}) < 1+ \lambda_0^2.
\eeq
\et
\begin{proof}
This estimate relies on the foliation structure of the {\af} manifold $M$ in an essential way. Our
strategy is the following: we construct a Fuchsian manifold $N$ from the {\ms} $\Sigma \subset M$
(assigning the warped product metric), and it is quasi-isometric to $M$. We then lift this
quasi-isomorphism to $\H^3$ and estimate the quasi-conformal constant in the cover.

Since the normal bundle over $\Sigma$ in $M$ is trivial, i.e., the geodesics perpendicular to $\Sigma$
are disjoint from each other. Therefore, any point $p\in{}M$ can be represented by the pair $p=(x,r)$, here
$x$ is the projection of $p$ to $\Sigma$ along the geodesic which passes through $p$ and is
perpendicular to $\Sigma$, and $r$ is the (signed) distance between $p$ and $x$. Now we can construct
a Fuchsian manifold $N=\Sigma \times\R$ as follows: suppose that the induced metric on
$\Sigma \subset{}M$ is given by $g(x)=e^{2v(x)}\I$, here $v(x)$ is a smooth function defined on $\Sigma$
and $\I$ is the $2\times{}2$ identity matrix. Let $\tilde{g}(x)$ be the (unique) {\hym} in the conformal class
of $g(x)$, then the (warped product) metric $\bar\rho$ on $N$ is given by
\begin{equation*}
   \bar\rho(x,r)=
   \begin{pmatrix}
      \cosh^2(r)\tilde{g}(x) & 0 \\
      0         & 1
   \end{pmatrix}\ ,
\end{equation*}
or $\bar\rho(x,r)= \cosh^2(r)\tilde{g}(x)+ dr^2$.
\medskip
Note that the surface $\Sigma\times\{0\}$ is totally geodesic. Similarly, any point $q\in{}N$ can be represented
by $q=(y,s)$, here $y$ is the projection of $q$ to $\Sigma\times\{0\}$ and $s$ is the distance between $q$
and $y$.

Now we may define a map $\varphi:N\to M$ by $\varphi(x,r)=(x,r)$ for $(x,r)\in N$. By the result in
\cite[p. 162]{Uhl83}, the map $\varphi$ is a quasi-isometry. We lift $\varphi$ to the map
$\tilde\varphi:\H^3\to\H^3$, then $\tilde\varphi$ is also a quasi-isometry. By the results in
\cite[Theorem 9]{Geh62}, \cite[Theorem 12.1]{Mos68}, \cite[Corollary 5.9.6]{Thu82} and
\cite[Theorem 3.22]{MT98}, the map $\tilde\varphi$ can be extend to an automorphism
\be
   \breve\varphi:\H^3\cup\widehat{\C}\to\H^3\cup\widehat{\C}
\ene
such that the restriction $\breve\varphi|{\widehat\C}=:f$ is a quasi-conformal mapping. In particular, $f$ maps
$S^{2}_{\pm}$ to $\Omega_{\pm}(\Gamma)$, where $S^{2}_{\pm}=S^{2}_{\infty}\setminus{}S^{1}$ are
hemispheres such that $\partial{}S^{2}_{+}=S^{1}=\partial{}S^{2}_{-}$, and $f(S^{1})=\Lambda_\Gamma$,
respectively.

We claim that $f|_{S^{2}_{+}}:S^{2}_{+}\to\Omega_{+}(\Gamma)$ is a $k$-quasiconformal mapping, with the
dilatation $K = \frac{1+k}{1-k}$, and
\be\label{K}
   K<\frac{1+\lambda_{0}}{1-\lambda_{0}}\ .
\ene
\medskip
To see this, we let $\Pi$ be the lift of the totally geodesic surface $\Sigma \times\{0\}\subset{}N$, and
$\widehat{\Sigma}$ be the lift of the surface $\Sigma\subset{}M$. Recall that the identity map between $\Pi$
and $\widehat{\Sigma}$ is an isometry, and we can define hyperbolic Gauss maps $G'_{+}:\Pi\to{}S^{2}_{+}$
and $G''_{+}:\widehat{\Sigma} \to\Omega_{+}(\Gamma)$ (as in \cite{Eps86}) such that we have the following
commutative diagram
\beq
\begin{CD}
   S^{2}_{+} @>{f}>> \Omega_{+}(\Gamma) \\
    @A{G'_{+}}AA    @AA{G''_{+}}A \\
   \Pi @>{\id}>> \widehat{\Sigma}
\end{CD}
\eeq
\vskip 0.1in
Since $G'_{+}$ is a conformal mapping, and $\id$ is an isometry, we therefore find that $G''_{+}\circ{}f$ is also
a conformal mapping. By Proposition 5.1 and Corollary 5.3 in \cite{Eps86}, $(G''_{+})^{-1}$ is a
$k$-quasiconformal mapping, and so is $f$. In particular, $\Lambda_{\Gamma}=f(S^{1})$ is a $k$-quasicircle.

Recently, Smirnov (\cite{Smi10}) proved Astala's conjecture (\cite{Ast94}): the {\Hd} of a $k$-quasicircle is at
most $1+k^2$. Now combining with \eqref{K}, we have
\begin{equation*}
k = \frac{K -1}{K+1} < \lambda_0.
\end{equation*}
Now our estimate follows easily.
\end{proof}


Note that this estimate partially answers a question raised by Uhlenbeck, see [Problem 5, Page 160, \cite{Uhl83}].

\section{Non-foliation for non-{\af} manifolds}\label{sec:counterexam}
Foliations of closed surfaces of {\cmc} play important role in three-dimensional geometry (\cite{Thu82})
and physics (for instance \cite{AMT97}). Recently Mazzeo-Pacard (\cite{MP07}) showed the existence of such a foliation
near either end of a {\qf} manifold. Note that, every AdS space-time admits such a foliation (\cite{BBZ07}), so a natural
question (given the analog from the AdS space-time) is to ask whether a {\qf} manifold admits a global such foliation. In
this section, we answer this question negatively. We construct a {\qf} manifold does not admit a foliation of closed
surfaces of {\cmc}, and we stress that this example is not {\af} since our construction admits at least two closed {\ms}s.

\subsection{Preparation}

The main scheme consists of the following three steps:
\begin{enumerate}
\item
we construct a {\qf} manifold $M$ obtained from $\H^3$ modulo a {\qf} group generated by reflections about some
circles on the Riemann sphere $S_{\infty}^{2}$ (see \S 4.2).
\item
For a given foliation of closed surfaces of {\cmc} on $M$, we lift it up to $\H^3$ where this foliation becomes a
foliation of hypersurfaces of {\cmc} which share the same asymptotic infinity. But for these circles on $S_{\infty}^{2}$,
one can construct cylinder-like {\ms}s which serve as barriers to force two leaves $L_{t_1}$ and $L_{t_2}$ of the
foliation of $\H^3$ become minimal (see \S 4.3).
\item
In the region of $\H^3$ bounded by $L_{t_1}$ and $L_{t_2}$, we find a leaf $L_{t_3}$ of small {\mc}
$H_0 = 2\tanh(\varepsilon)$. Then we use two disks $D_{1}(\varepsilon)$ and $D_{2}(\varepsilon)$, both of {\cmc} $H_0$, to push
the leaf $L_{t_3}$ to self-intersect.
\end{enumerate}

We will work in the ball model of $\H^{3}$, i.e.,
\beq
   \H^{3}=\{(x,y,z)\in\R^{3}\ |\
   x^{2}+y^{2}+z^{2}<{}1\},
\eeq
\medskip
equipped with metric
\beq
   ds^{2}=\frac{4(dx^{2}+dy^{2}+dz^{2})}{(1-r^{2})^{2}}\ ,
\eeq
\medskip
where $r=\sqrt{x^{2}+y^{2}+z^{2}}$.

The hyperbolic space $\H^{3}$ has a natural compactification: $\overline{\H}{}^{3}=\H^{3}\cup{}S_{\infty}^{2}$,
where $S_{\infty}^{2}=\widehat{\C}$ is the Riemann sphere. Suppose $X$ is a subset of $\H^{3}$, we define
$\partial_{\infty}X$ by
\beq
   \partial_{\infty}X=\overline{X}
   \cap{}S_{\infty}^{2}\ ,
\eeq
the {\em asymptotic boundary} of $X$, where $\overline{X}$ is the closure of $X$ in $\overline{\H}{}^{3}$.

In anticipation of the barrier surfaces that we will use later, we need some results of Gomes and L{\'o}pez
(see \cite{Gom87,Lop00}). Let us first make some definitions:

Suppose $G$ is a subgroup $\Isom(\H^{3})$ which leaves a geodesic $\gamma\subset\H^{3}$ pointwisely fixed.
We call $G$ the {\it spherical group} of $\H^{3}$ and $\gamma$ the {\it rotation axis} of $G$. A surface in
$\H^{3}$ invariant under $G$ is called a {\em spherical surface}. For two circles $C_{1}$ and $C_{2}$ in
$\H^{3}$, if there is a geodesic $\gamma$, such that each of $C_{1}$ and $C_{2}$ is invariant under the group
of rotations that fixes $\gamma$ pointwisely, then $C_{1}$ and $C_{2}$ are said to be {\em coaxial}, and
$\gamma$ is called the {\em rotation axis} of $C_{1}$ and $C_{2}$.

Let $P_{1}$ and $P_{2}$ be two disjoint geodesic plane in $\H^{3}$, then $P_{1}\cup{}P_{2}$ divides $\H^{3}$
into three components. Let $X_{1}$ and $X_{2}$ be the two of them with $\partial{}X_{i}=P_{i}$ for $i=1,2$.
Given two subsets $A_{1}$ and $A_{2}$ of $\overline{\H}{}^{3}$, we say $P_{1}$ and  $P_{2}$ {\em separate}
$A_{1}$ and $A_{2}$ if one of the following cases occurs (\cite{Lop00}):
\begin{enumerate}
   \item if $A_{1},A_{2}\subset\H^{3}$, then
         $A_{i}\subset{}X_{i}$ for $i=1,2$;
   \item if $A_{1}\subset\H^{3}$ and
         $A_{2}\subset{}S_{\infty}^{2}$, then
         $A_{1}\subset{}X_{1}$ and
         $A_{2}\subset\partial_{\infty}X_{2}$;
   \item if $A_{1},A_{2}\subset{}S_{\infty}^{2}$, then
         $A_{i}\subset\partial_{\infty}X_{i}$ for $i=1,2$.
\end{enumerate}
\medskip
Then we may define the distance between $A_{1}$ and $A_{2}$ by
\be\label{eq:distance between circles}
   d(A_{1},A_{2})=\sup\{\dist(P_{1},P_{2})\ |\
   P_{1}\ \text{and}\ P_{2}\ \text{separate}\ A_{1}\
   \text{and}\ A_{2}\}\ ,
\ene
where $\dist(P_{1},P_{2})$ is the hyperbolic distance between
$P_{1}$ and $P_{2}$. We need the following result
of Gomes to ensure the existence of a {\ms} with
$C_{1}\cup{}C_{2}$ as its asymptotic boundary. Namely,

\bl[\cite{Gom87}]
\label{lem:Gomes1987}
There exists a finite constant $d_{0}>0$ such that
for two disjoint circles $C_{1},C_{2}\subset{}S_{\infty}^{2}$, if
$d(C_{1},C_{2})\leq{}d_{0}$, then there exists a {\ms}
$\Pi$ which is a surface of revolution with asymptotic
boundary $C_{1}\cup{}C_{2}$.
\el

We will call the minimal surface $\Pi$ in Lemma~\ref{lem:Gomes1987}
a {\em minimal catenoid}.

Next we need a similar result for surfaces of {\cmc}.
To do this, we let $C_{1}$ and $C_{2}$ be two disjoint circles
on $S_{\infty}^{2}$, and let $P_{1}$ and $P_{2}$ be two
geodesic planes whose asymptotic boundaries are
$C_{1}$ and $C_{2}$, respectively. Suppose $C'_{1}\subset{}P_{1}$
and $C'_{2}\subset{}P_{2}$ such that $C'_{1}$
and $C'_{2}$ are two coaxial circles with respect to the rotation
axis of $C_{1}$ and $C_{2}$. The following result
is due to L\'{o}pez:

\bl[\cite{Lop00}]
\label{lem:Lopez2000}
Given a constant $H\in(-2,2)$, there exists a constant $d_{H}$,
depending only on
$H$ such that if $d(C_{1},C_{2})\leq{}d_{H}$,
then there exists a compact smooth surface $\Pi'$ such that
$\partial\Pi'=C'_{1}\cup{}C'_{2}$
and the mean curvature of $\Pi'$ is equal to $H$ with
respect to the inward normal vector, i.e., the normal vector
pointing to the domain containing the rotation axis of $C_{1}$ and $C_{2}$.
\el

\begin{rem}
In Lemma~\ref{lem:Lopez2000}, when $H<0$, then there is no such a surface $\Pi'$ of {\cmc}
if we replace $C'_{i}$ by $C_{i}$ for $i=1,2$ (see \cite{Pal99}). For the compact surface $\Pi'$ in
Lemma~\ref{lem:Lopez2000}, it might not be a surface of revolution according to the discussion
in \cite[p. 234]{Lop00}. But it is sufficient for our goals.
\end{rem}

\subsection{The Construction} 
In this subsection, we complete the first step mentioned in \S 4.1, namely, we construct a {\qf} manifold by using groups
of reflections about circles on $S_{\infty}^{2}$.

\begin{figure}[ht]
\centering
      \includegraphics[scale=0.6]{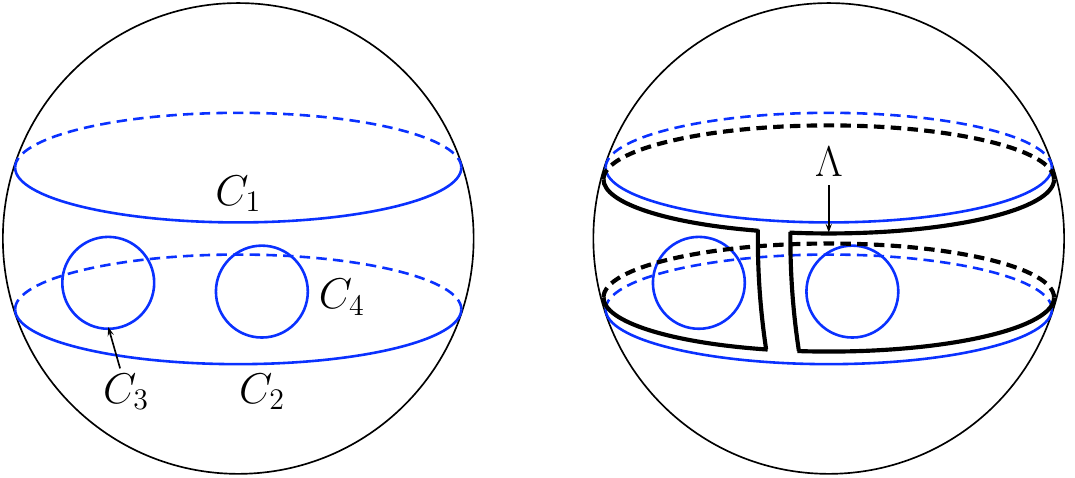}
      \caption{Four Circles and a ``narrow bridge"}\label{fig: four circles}
\end{figure}

Let $\varepsilon>0$ be a sufficiently small number such that
$\varepsilon\ll{}d_{0}/2$, here $d_{0}$ is the constant in
Lemma~\ref{lem:Gomes1987}, and let
\begin{equation}
   H_{0}=2\tanh\varepsilon\ .
\end{equation}
Consider $\H^3$ as a unit ball in $\R^3$ and
consider $S_{\infty}^{2}$ as a unit sphere in $\R^3$.
We pick up four circles $\{C_{i}\}_{i=1,\ldots,4}$ on
$S_{\infty}^{2}$ and four geodesic planes $\{D_{i}\}_{i=1,\ldots,4}$
in $\H^3$ as follows (see Figure~\ref{fig: four circles}).
\begin{enumerate}
   \item Let $C_1$ and $C_2$ be the circles on the
         horizontal planes $z=\tanh\varepsilon$ and
         $z=-\tanh\varepsilon$ respectively. It's easy to verify that
         $d(C_{1},C_{2})=2\varepsilon$,
         where $d(\cdot,\cdot)$ is the distance defined
         by \eqref{eq:distance between circles}.
   \item Let $C_{3}$ and $C_{4}$ be two disjoint circles between the horizontal
         planes $z=\tanh\varepsilon$ and $z=-\tanh\varepsilon$ such that
         \begin{itemize}
            \item $C_3$ and $C_4$ have the same size with respect the
                  spherical metric on the asymptotic boundary $S_{\infty}^{2}$,
                  and
            \item the distance between $C_3$ and $C_4$ with respect the
                  spherical metric is sufficiently small so that
                  $d(C_{3},C_{4})\ll{}\min\{d_{H_{0}},d_{0}\}$.
         \end{itemize}
   \item Let $D_{i}$ be the geodesic plane in $\H^{3}$ that is asymptotic
         to $C_i$, i.e., $\partial_{\infty}D_{i}=C_{i}$ for $i=1,\ldots,4$.
\end{enumerate}
By the above construction, $d(C_1,C_2)=\dist(D_1,D_2)$ and
$d(C_3,C_4)=\dist(D_3,D_4)$.

Let $\Lambda$ be a closed piecewise smooth curve on $S_{\infty}^{2}$ which
disconnect $C_1\cup{}C_2$ from $C_3\cup{}C_4$
(see Figure~\ref{fig: four circles}), then we cover $\Lambda$ by finitely many
disks $\{B_{l}\subset{}S_{\infty}^{2}\}_{l=1,\ldots,N}$ with
small radii such that
\ben
   \item each circle $\partial{}B_{l}$ is invariant under the
         rotation along the geodesic connecting the
         origin $O$ and the center of the disk $B_{l}$, which locates at $\Lambda$;
   \item the radii of disks are small enough such that
         $B_{l}\cap{}C_{i}=\emptyset$ for $l=1,\ldots,N$ and
         $i=1,\ldots,4$; and
   \item for each $l\equiv{}1\ (\mod{}N)$, $\partial{}B_{l}$ intersects
         both $\partial{}B_{l-1}$ and $\partial{}B_{l+1}$ perpendicularly
         and no other circles,
\een
then we obtain a torsion free group $\Gamma$ which is the subgroup of
orientation preserving transformations in the group
generated by $N$ reflections about the circles
$\partial{}B_{1},\ldots,\partial{}B_{N}$. It is well-known that such
a subgroup is {\qf} (\cite [Page 263]{Ber72}).
Therefore the {\htm} $M=\H^3/\Gamma$ is a {\qf} manifold. The limit
set of the {\qf} group $\Gamma$, denoted by $\Lambda_{\Gamma}$,
is around the curve $\Lambda$. Let
$S_{\infty}^{2}\setminus\Lambda_{\Gamma}=\Omega_{\pm}$, where $\Omega_{-}$
contains $C_{1}\cup{}C_{2}$ and $\Omega_{+}$ contains
$C_{3}\cup{}C_{4}$ (See Figure~\ref{fig: four circles}).
\vskip 0.1in
\subsection{The Proof}
We will first recall the Hopf's maximum principle for tangential hypersurfaces in Riemannian geometry,
which will be used in the proof of the Theorem \ref{counter}. We will actually
only compare hypersurfaces in $\H^3$ which are of {\cmc}.

\bl [\cite{Hop89}]\label{Hopf}
Let $S_{1}$ and $S_{2}$ be two hypersurfaces in a Riemannian
manifold which intersect at a common point
tangentially. If $S_{2}$ lies in positive side of $S_{1}$ around
the common point, then $H_{1} \le H_{2}$, where
$H_{i}$ is the {\mc} of $S_{i}$ at the common point for $i=1,2$.
\el

For our particular situation, we need a corollary of the maximum principle that 
will be used later. 

\begin{cor}\label{cor:disjoint CMC surfaces}
Let $\Sigma\subset\H^3$ be a disk-type surface whose asymptotic
boundary is a Jordan curve and its mean curvature is a constant $H\in(-2,2)$. Let
$D$ be a totally geodesic plane in $\H^3$, and let $D(r)$ be one of the components
of the boundary of $\Nscr_{r}(D)$ for $r=\tanh^{-1}(|H|/2)$, then $D(r)$ is a disk-type surface with the same
{\cmc} as that of $\Sigma$, here $\Nscr_{r}(D)$ denotes the 
$r$-neighborhood $D$. Suppose that the normal vectors on
$\Sigma$ and $D(r)$ are in the same direction.
If $\partial_{\infty}D\cap\partial_{\infty}\Sigma=\emptyset$, then
$D(r)\cap\Sigma=\emptyset$.
\end{cor}

\begin{proof}Let $\partial_{\infty}D=C$ and $\partial_{\infty}\Sigma=\Lambda$.
By a M\"obius transformation, we may assume that $C$ is on the horizontal
plane $z=0$, i.e. $C=\{(x,y,0)\in\R^3\ |\ x^2+y^2=1\}$ and that $\Lambda$ is above
the the horizontal plane $z=0$. Besides, we also assume that the normal vectors on
$D(r)$ and $\Sigma$ are {\em downward}, i.e. the normal vectors point to the domains which
are below the surfaces $D(r)$ and $\Sigma$ respectively.

Let $W$ be the subdomain of $\H^3$ which is below $D(r)$. Note that the surface $D(r)$ has constant {\pc}. 
Using the translations along the $z$-axis, we may foliate $W$
by disk-type surfaces $\{D_{t}(r)\}_{-\infty<t<0}$ whose asymptotic
boundaries are circles and whose mean curvatures are the same as that of $D(r)$
with respect to the {\em downward} normal vectors.
If $D(r)\cap\Sigma\ne\emptyset$, then some interior points of $\Sigma$ are contained
in $W$, so from bottom to top, there is a surface $D_{t}$ touches $\Sigma$ for the
first time, here $t\in(-\infty,0)$, therefore $H(D_t(r))>H(\Sigma)$ by the 
maximal principle, here $H(\cdot)$ denotes
the mean curvature of the surface with respect to the downward
normal vectors. This is impossible, since they are supposed to be equal. Thus
$D(r)$ must be disjoint from $\Sigma$.
\end{proof}

Now we are ready to show:

\bt\label{counter}
The {\qf} manifold $M=\H^{3}/\Gamma$ constructed above can not be
foliated by closed surfaces of {\cmc}.
\et

\bp
We will argue by contradiction and we follow the scheme outlined in
the beginning of \S 4.1. Let us assume
that our construction $M$ is foliated by surfaces of {\cmc},
where each surface is closed and incompressible.
We lift this foliation to the universal covering space
$\H^{3}$, then we obtain a foliation of $\H^{3}$ such that each leaf
is a disk of {\cmc} and all disks share the same
asymptotic boundary $\Lambda_{\Gamma}$.
\vskip 0.1in
\noindent
Step $(2)$: {\it Existence of two minimal leaves $L_{t_1}$ and $L_{t_2}$}.
\vskip 0.1in

Notice that any disk-type surface in $\H^{3}$ with asymptotic boundary
$\Lambda_{\Gamma}$ divides $\overline{\H}{}^{3}$ into two parts, one of them contains $\Omega_{-}$,
the other contains $\Omega_{+}$. We choose the normal vector field on the
disk-type surface so that each normal vector points to the domain
containing $\Omega_{-}$. Assume that there is a {\cmc}
foliation $\Fcal=\{L_{t}\}$ of $\H^3$, with a parameter
$t\in(-\infty,\infty)$ such that the leaves are convergent to
$\Omega_{\pm}$ as $t\to{}\pm\infty$ respectively.
In other words, we have
\be\label{eq:foliation limit behavior}
   \lim_{t\to\pm\infty}H(L_{t})=\pm{}2\ ,
\ene
where $H(L_{t})$ denotes the {\mc} of the leaf $L_{t}$ with
respect to the normal vectors pointing to the domain
containing $\Omega_{-}$. Here we just need to assume that
$H(L_{t})$ is a continuous function
of the parameter $t\in(-\infty,\infty)$.

\begin{figure}[ht]
      \includegraphics[scale=0.6]{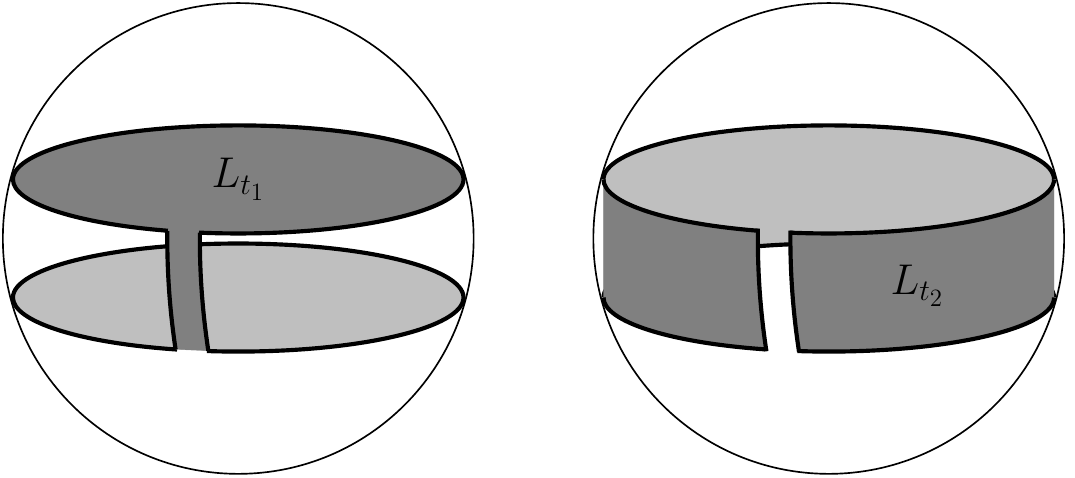}
      \caption{Two Minimal Disks}\label{fig: minimal surfaces}
\end{figure}

Since $d(C_{3},C_{4})\ll{}\min\{d_{H_{0}},d_{0}\}$,
there exists a minimal catenoid $\Pi_{1}$
whose asymptotic boundary is $C_{3}\cup{}C_{4}$ by
Lemma~\ref{lem:Gomes1987}. Starting from $\Omega_{-}$ to $\Omega_{+}$,
there is a leaf $L_{t'}\in\Fcal$
which touches $\Pi_{1}$ for the first time,
then the {\mc} of the leaf $L_{t'}$ must be positive by the
maximal principle (Lemma \ref{Hopf}). Because of the
limiting behavior in \eqref{eq:foliation limit behavior},
there exists $t_1\in(-\infty,t')$ such that the {\mc} of $L_{t_{1}}$
is zero, i.e. the leaf $L_{t_{1}}$ is a disk type {\ms}
(see the left figure in Figure~\ref{fig: minimal surfaces}).
Besides, we may choose $t_1$ small enough such that
$H(L_t)<0$ for all $t\in(-\infty,t_1)$; this can be done since
all leaves in $\Fcal$ which are minimal must be contained
in the convex hull of $\Lambda_{\Gamma}$.

Similarly, since $d(C_{1},C_{2})\ll{}d_{0}$, there exists a
minimal catenoid $\Pi_{2}$ whose asymptotic boundary is
$C_{1}\cup{}C_{2}$ by
Lemma~\ref{lem:Gomes1987}. Starting from $\Omega_{+}$ to $\Omega_{-}$,
there exists a leaf $L_{t''}\in\Fcal$
which touches $\Pi_{2}$ for the first time,
then the {\mc} of the leaf $L_{t''}$ must be negative
(Lemma \ref{Hopf} again). Because of the
limiting behavior in \eqref{eq:foliation limit behavior},
there exists $t_{2}\in(t'',\infty)$ such that the {\mc} of $L_{t_{2}}$
is zero, i.e. the leaf $L_{t_{2}}$ is a disk type {\ms}
(see the right figure in Figure~\ref{fig: minimal surfaces}).
We may choose $t_2$ big enough such that
$H(L_t)>0$ for all $t\in(t_2,\infty)$.
\vskip 0.1in

By the construction, $L_{t_{1}}\ne{}L_{t_{2}}$. Besides, $L_{t_{1}}$
is close to $\Omega_{-}$ and $L_{t_{2}}$ is close to $\Omega_{+}$,
so we must have $t_1<t_2$. Note that we obtain two distinct closed
{\ms}s in the quotient manifold $M$, and hence $M$ can not be {\af}.
\vskip 0.1in
\noindent
Step $(3)$: {\it Some leaf $L_{t_3}$ must self-intersect}.
\vskip 0.1in
Let $X\subset\H^{3}$ be the region bounded by leaves $L_{t_{1}}$ and
$L_{t_{2}}$, then by assumption $X$ is
foliated by $\{L_{t}\}_{t_{1}\leq{}t\leq{}t_{2}}$, i.e.
\beq
   X=\bigcup_{t_{1}\leq{}t\leq{}t_{2}}L_{t}\ .
\eeq
Let
\beq
   W_{1}=\bigcup_{-\infty<t\leq{}t_{1}}L_{t}
   \qquad\text{and}\qquad
   W_{2}=\bigcup_{t_{2}\leq{}t<\infty}L_{t}\ .
\eeq
Clearly we have $\H^3=W_1\cup{}X\cup{}W_2$, $W_1\cap{}W_2=\emptyset$
and $X\cap{}W_{i}=L_{t_i}$ for $i=1,2$.

Recall that both $D_{3}$ and $D_{4}$ are totally geodesic planes, and
$\partial_{\infty}D_{i}\cap\Lambda_{\Gamma}=
C_{i}\cap\Lambda_{\Gamma}=\emptyset$
for $i=3,4$, so both $D_{3}$ and
$D_{4}$ are disjoint from $L_{t_{2}}$
by Corollary \ref{cor:disjoint CMC surfaces}
(In fact, $D_{i}\subset{}W_{2}$ for $i=3,4$).
We choose two circles $C_{3}'\subset{}D_{3}$ and
$C_{4}'\subset{}D_{4}$ such that $C_{3}'$ and $C_{4}'$ are
coaxial with respect to the rotation axis of $C_{3}$
and $C_{4}$. By the Lemma~\ref{lem:Lopez2000}, there is a compact
surface $\Pi_{0}$ of {\cmc} $-H_{0}$ with respect
to the {\em inward normal vectors}, i.e., the normal vectors
pointing to the domain containing the rotation axis of $C_{3}'$
and $C_{4}'$.

We claim that $\Pi_{0}\cap{}W_{1}=\emptyset$. Since
$D_{i}\subset{}W_{2}$ for $i=3,4$, then
$C_{i}'\subset{}W_{2}$ for $i=3,4$. Therefore, if
$\Pi_{0}\cap{}W_{1}\ne\emptyset$, then any point
in $\Pi_{0}\cap{}W_{1}$ must be the interior point of
$\Pi_{0}$. Starting from $\Omega_{-}$ to $L_{t_1}$, let $L_{t}$ be
the leaf contained in $W_{1}$
which touches $\Pi_{0}$ for the first time, then
$H(L_{t})>H_{0}$, here the normal vector at the common point points
to the domain containing $\Omega_{-}$, which is the {\em outward} normal
vector on $\Pi_{0}$. This is impossible since $H(L_{t})\leq{}0$ for
all $t\in(-\infty,t_1]$.

In particular, we have $L_{t_1}\cap\Pi_0=\emptyset$ and
$L_{t_2}\cap\Pi_0\ne\emptyset$. Therefore, starting from
$L_{t_1}$ to $L_{t_2}$, there exists $t_{*}\in(t_1,t_2)$ such that
the leaf $L_{t_*}$ touches $\Pi_{0}$
for the first time, then $H(L_{t_*})>H_{0}$
by the maximal principle, here the normal vector at the common point
points to the domain containing $\Omega_{-}$. So there exists
$t_{3}\in(t_{1},t_{*})$ such that $H(L_{t_{3}})=H_{0}$. Notice
that $L_{t_{3}}$ is {\em close} to $L_{t_{1}}$, so its shape is similar to
that of $L_{t_{1}}$, one may imagine that $L_{t_{3}}$ consists of two
disks connected by a narrow bridge. Besides, $L_{t_{3}}$ is still contained
in $X$.

We now complete the step $(3)$, namely, the leaf $L_{t_{3}}$ must self-intersect.
To see this, let $D_{1}(\varepsilon)$ be the component of
$\partial\Nscr_{\varepsilon}(D_1)$ that is below the geodesic
plane $D_1$, here $\Nscr_{\varepsilon}(D_1)$ is the
(hyperbolic) $\varepsilon$-neighborhood of $D_{1}$,
then $D_{1}(\varepsilon)$ is a surface with
constant mean curvature $H_{0}$ with respect to the {\em upward}
normal vectors, i.e. the normal vectors pointing to domain not
containing $C_{2}$. It is well known that the equidistant surface from a totally geodesic
disk in $\H^3$ is of {\cmc}. Similarly, let $D_{2}(\varepsilon)$ be the
component of $\partial\Nscr_{\varepsilon}(D_2)$ that is above the
geodesic plane $D_2$, then $D_{2}(\varepsilon)$ is a surface with
constant mean curvature $H_{0}$ with respect to the {\em downward}
normal vectors, i.e. the normal vectors pointing to domain not
containing $C_{1}$. One can see $D_1(\varepsilon)$ as a dome below $D_1$, while
$D_2(\varepsilon)$ as a dome above $D_2$.

Since $\partial_{\infty}D_{i}\cap\partial_{\infty}L_{t_{3}}=
C_{i}\cap\Lambda_{\Gamma}=\emptyset$
for $i=1,2$, neither $D_{1}(\varepsilon)$ nor
$D_{2}(\varepsilon)$ intersects $L_{t_{3}}$ by
Corollary \ref{cor:disjoint CMC surfaces}.
Recall that the shape of $L_{t_{3}}$ is similar to that of $L_{t_{1}}$,
i.e., $L_{t_{3}}$ consists of two disks connected by a narrow bridge.
The two disks of $L_{t_{3}}$ must be below $D_{1}(\varepsilon)$ and
above $D_{2}(\varepsilon)$.
Recall that the (hyperbolic) distance between $D_1$ and
$D_2$ is $2\varepsilon$, therefore
$D_{1}(\varepsilon)\cap{}D_{2}(\varepsilon)=\{O\}$, where $O\in\H^{3}$
is the origin, so $L_{t_{3}}$ must self intersect.
Now the theorem follows easily.
\ep

\bibliographystyle{amsalpha}
\bibliography{ref-af}
\end{document}